\documentclass{gtart_h}

\def\ifplaintex{\expandafter\ifx\csname documentclass\endcsname\relax}

\def\gtp{{\mathsurround=0pt\it $\cal G\mskip-2mu$eometry \&\ 
$\cal T\!\!$opology $\cal P\!$ublications}}  

\def\recd{{\small Received:\qua\receiveddate\ifx\reviseddate\relax
\else\qquad Revised:\qua\reviseddate\fi\par}} 

\def\lognumber#1{\def\thelognumber{#1}}
\def\volumenumber#1{\def\thevolumenumber{#1}}
\def\volumeyear#1{\def\thevolumeyear{#1}}
\def\papernumber#1{\def\thepapernumber{#1}}
\def\pagenumbers#1#2{\def\startpage{#1}\def\finishpage{#2}}
\def\published#1{\def\publishdate{#1}}

\def\received#1{\def\receiveddate{#1}}

\def\accepted#1{\def\accepteddate{#1}}
\def\asciititle#1{\def\theasciititle{#1}}

\def\asciiauthors#1{\def\theasciiauthors{#1}}
\def\asciiaddress#1{\def\theasciiaddress{#1}}
\def\asciiemail#1{\def\theasciiemail{#1}}

\def\coverauthors#1{\def\thecoverauthors{#1}}
\long\def\asciiabstract#1{\long\def\theasciiabstract{#1}}

\let\\\par\let\thelognumber\relax\let\thevolumenumber\relax
\let\thepapernumber\relax\let\thevolumeyear\relax\let\startpage\relax
\let\finishpage\relax\let\publishdate\relax\let\receiveddate\relax
\let\reviseddate\relax\let\accepteddate\relax\let\theasciititle\relax
\let\theasciiauthors\relax\let\theasciiaddress\relax
\let\theasciiabstract\relax

\let\thecoverauthors\relax\let\theasciiemail\relax

\ifplaintex
\font\logobig=cmssbx10 scaled 3836
\font\logomed=cmssbx10 scaled 2557
\else
\font\logobig=cmssbx10 scaled 4200
\font\logomed=cmssbx10 scaled 2800
\fi

\long\def\makeagttitle{   
\count0=\startpage
\agt\hfill      
\hbox to 45truept{\vbox to 0pt{\vglue -13truept{\logomed A\kern -.37em{\logobig 
T}\kern -.38em G}\vss}\hss}
\break
{\small Volume \thevolumenumber\ (\thevolumeyear)
\startpage--\finishpage\nl
Published: \publishdate}

\vglue .25truein

{\parskip=0pt\leftskip 0pt plus
1fil\def\\{\par\smallskip}{\Large\bf\thetitle}\par\medskip} \vglue
0.05truein

{\parskip=0pt\leftskip 0pt plus 1fil\def\\{\par}{\sc\theauthors}
\par\medskip}%
 
\vglue 0.03truein 

{\small\leftskip 25truept\rightskip 25truept{\bf Abstract}\stdspace\theabstract

{\bf AMS Classification}\stdspace\theprimaryclass
\ifx\thesecondaryclass\relax\else; \thesecondaryclass\fi\par
{\bf Keywords}\stdspace \thekeywords\par}\vglue 7truept

}   

\ifplaintex
\font\phead=cmsl9 scaled 950
\font\pnum=cmbx10 scaled 913
\font\pfoot=cmsl9 scaled 950
\headline{\vbox to 0pt{\vskip -4.5mm\line{\small\phead\ifnum
\count0=\startpage ISSN 1472-2739 (on-line) 1472-2747 (printed)
\hfill {\pnum\folio}\else\ifodd\count0\def\\{ }%
\ifx\theshorttitle\relax\thetitle\else\theshorttitle\fi\hfill{\pnum\folio}
\else\def\\{ and }{\pnum\folio}\hfill\ifx\theshortauthors\relax\theauthors
\else\theshortauthors\fi\fi\fi}\vss}}
\footline{\vbox to 0pt{\vglue 0mm\line{\small\pfoot\ifnum\count0=\startpage
\copyright\ \gtp\hfill\else
\agt, Volume \thevolumenumber\ (\thevolumeyear)\hfill\fi}\vss}}
\else
\font=cmsl9 scaled 1050
\font\lnum=cmbx10 
\font=cmsl9 scaled 1050
\makeatletter
\def\@oddhead{{\small\lhead\ifnum\count0=\startpage ISSN 1472-2739 
(on-line) 1472-2747 (printed)\hfill {\lnum\number\count0}\else\ifodd\count0
\def\\{ }\ifx\theshorttitle\relax \thetitle \else\theshorttitle\fi\hfill
{\lnum\number\count0}\else\def\\{ and }{\lnum\number\count0}
\hfill\ifx\theshortauthors\relax 
\theauthors\else\theshortauthors\fi\fi\fi}}\def\@evenhead{\@oddhead}
\def\@oddfoot{\small\lfoot\ifnum\count0=\startpage\copyright\ \gtp\hfill\else
\agt, Volume \thevolumenumber\ (\thevolumeyear)\hfill\fi}
\def\@evenfoot{\@oddfoot}
\makeatother
\fi
\let\maketitlepage\makeagttitle
\let\makeshorttitle\maketitlepage
\let\maketitle\maketitlepage

\newwrite\gtoutfile
\long\gdef\makeheadfile{  
{\def\\{, }\def\s{ }
\immediate\openout\gtoutfile head.xxx
\immediate\write\gtoutfile{Proxy-for: \ifx\theasciiauthors\relax
\theauthors\else\theasciiauthors\fi\s<\ifx\theasciiemail\relax\theemail\else\theasciiemail\fi>}
\immediate\write\gtoutfile{\noexpand\\}
\immediate\write\gtoutfile{Authors: \ifx\theasciiauthors\relax
\theauthors\else\theasciiauthors\fi}
{\def\\{ }\immediate\write\gtoutfile{Title: \ifx\theasciititle\relax
\thetitle\else\theasciititle\fi}}
\immediate\write\gtoutfile{Subj-class: GT or SG, GR etc}
\immediate\write\gtoutfile{MSC-class: \theprimaryclass\ifx\thesecondaryclass\relax\else, \thesecondaryclass\fi}
\immediate\write\gtoutfile{Journal-ref: Algebr. Geom. Topol. \thevolumenumber\s
(\thevolumeyear) \startpage-\finishpage}
\immediate\write\gtoutfile{Comments: Published by Algebraic and
Geometric Topology at}
\immediate\write\gtoutfile{\s\s\s  http://www.maths.warwick.ac.uk/agt/AGTVol\thevolumenumber/agt-\thevolumenumber-\thepapernumber.abs.html}
\immediate\write\gtoutfile{\noexpand\\}
\immediate\write\gtoutfile{}
\ifx\theasciiabstract\relax
\immediate\write\gtoutfile{\theabstract}\else
\immediate\write\gtoutfile{\theasciiabstract}\fi
\immediate\write\gtoutfile{}
\immediate\write\gtoutfile{\noexpand\\}
\immediate\write\gtoutfile{}
\immediate\closeout\gtoutfile}}  

\def\maketitlepage{\makeagttitle\makeheadfile}
\let\makeshorttitle\maketitlepage
\let\maketitle\maketitlepage

\lognumber{29}
\volumenumber{5}
\volumeyear{2005}
\papernumber{29}
\pagenumbers{713}{724}
\received{13 February 2005} 
\accepted{30 June 2005}
\published{5 July 2005}

\usepackage[all]{xy}
\usepackage{amsmath,amssymb}
\theoremstyle{plain}
\newtheorem{theorem}{Theorem}
\newtheorem{proposition}[theorem]{Proposition}

\theoremstyle{definition}
\newtheorem{definition}[theorem]{Definition}
\newtheorem{example}[theorem]{Example}

\newcommand{\secref}[1]{Section~\ref{#1}}

\newcommand{\thmref}[1]{Theorem~\ref{#1}}
\newcommand{\propref}[1]{Proposition~\ref{#1}}

\newcommand{\defref}[1]{Definition~\ref{#1}}

\def\C{\mathbb C}
\def\Z{\mathbb Z}

\def\cF{\mathcal F}

\def\cat{{\rm cat}}
\def\wcatg{{\rm wcat_G}}
\def\clup{{\rm cup}}
\def\dim{{\rm dim}}
\def\conn{{\rm conn}}
\def\dl{{\rm dl}}

\begin{document}
\title{$H$-space structure on pointed mapping spaces}
\asciititle{H-space structure on pointed mapping spaces}

\author{Yves F\'elix\\Daniel Tanr\'e}
\asciiauthors{Yves Felix\\Daniel Tanre}
\coverauthors{Yves F\noexpand\'elix\\Daniel Tanr\noexpand\'e}

\address{D\'epartement de Math{\'e}matiques, Universit\'e Catholique de 
Louvain\\2, Chemin du Cyclotron, 1348 Louvain-La-Neuve, Belgium}

\secondaddress{D\'epartement de Math{\'e}matiques, UMR 8524,
Universit\'e de Lille~1\\59655 Villeneuve d'Ascq Cedex, France}

\asciiaddress{Departement de Mathematiques, Universite Catholique de 
Louvain\\2, Chemin du Cyclotron, 1348 Louvain-La-Neuve, 
Belgium\\and\\Departement de Mathematiques, UMR 8524,
Universite de Lille 1\\59655 Villeneuve d'Ascq Cedex, France}

\asciiemail{felix@math.ucl.ac.be, Daniel.Tanre@univ-lille1.fr}
\gtemail{\mailto{felix@math.ucl.ac.be}, \mailto{Daniel.Tanre@univ-lille1.fr}}

\begin{abstract} 
We investigate the existence of an $H$-space structure on the function
space, $\cF_*(X,Y,*)$, of based maps in the component of the trivial
map between two pointed connected CW-complexes $X$ and $Y$. For that,
we introduce the notion of $H(n)$-space and prove that we have an
$H$-space structure on $\cF_*(X,Y,*)$ if $Y$ is an $H(n)$-space and
$X$ is of Lusternik-Schnirelmann category less than or equal to
$n$. When we consider the rational homotopy type of nilpotent finite
type CW-complexes, the existence of an $H(n)$-space structure can be
easily detected on the minimal model and coincides with the
differential length considered by Y. Kotani.  When $X$ is finite,
using the Haefliger model for function spaces, we can prove that the
rational cohomology of $\cF_*(X,Y,*)$ is free commutative if the
rational cup length of $X$ is strictly less than the differential
length of $Y$, generalizing a recent result of Y. Kotani.
\end{abstract}

\asciiabstract{%
We investigate the existence of an H-space structure on the function
space, F_*(X,Y,*), of based maps in the component of the trivial
map between two pointed connected CW-complexes X and Y. For that,
we introduce the notion of H(n)-space and prove that we have an
H-space structure on F_*(X,Y,*) if Y is an H(n)-space and
X is of Lusternik-Schnirelmann category less than or equal to
n. When we consider the rational homotopy type of nilpotent finite
type CW-complexes, the existence of an H(n)-space structure can be
easily detected on the minimal model and coincides with the
differential length considered by Y. Kotani.  When X is finite,
using the Haefliger model for function spaces, we can prove that the
rational cohomology of F_*(X,Y,*) is free commutative if the
rational cup length of X is strictly less than the differential
length of Y, generalizing a recent result of Y. Kotani.}

\primaryclass{55R80, 55P62, 55T99}

\keywords{Mapping spaces, Haefliger model, Lusternik-Schnirelmann category}

\maketitle

\section{Introduction}
Let $X$ and $Y$ be pointed connected  CW-complexes. We study the occurrence of an $H$-space structure on the function space, $\cF_*(X,Y,*)$, of  based maps in the component of the trivial map. Of course when $X$ is a co-$H$-space or $Y$ is an $H$-space
this mapping space is an $H$-space. Here, we are considering weaker conditions, both on $X$ and $Y$, which guarantee the existence of an $H$-space structure on the function space. 
In \defref{def:nH}, we  introduce the notion of $H(n)$-space designed for this purpose and prove:

\begin{proposition}\label{prop:cat}
Let $Y$ be an $H(n)$-space and $X$ be a space of Lusternik-Schnirelmann category less than or equal to $n$. Then the space ${\mathcal F}_*(X,Y,*)$ is an $H$-space.
\end{proposition}

The existence of an $H(n)$-structure and the Lusternik-Schnirelmann category (LS-category in short) are hard to determine. We first study some properties of $H(n)$-spaces and give some examples. Concerning the second hypothesis, we are interested in replacing  $\cat(X)\leq n$ by an upper bound on an approximation of the LS-category (see \cite[Chapter~2]{CLOT03}).
We succeed in \propref{prop:wcatg} with an hypothesis on the dimension of $X$ but the most interesting replacement is obtained in the rational setting which constitutes the second part of this paper.

We use   Sullivan minimal models for which we refer to \cite{F-H-T99}. We recall here that each finite
type nilpotent CW-complex $X$ has a unique minimal model $(\land
V,d)$ that characterises all the rational homotopy type of $X$. 
We first prove that the existence of an $H(n)$-structure on a rational space $X_0$ can be easily detected from its minimal model. It corresponds to a valuation of the differential of this model,  introduced by Y.~Kotani in \cite{Kot04}:\\
The differential $d$ of the minimal model $(\land V,d)$ can be written as $d = d_1+d_2 + \cdots $
where $d_i$ increases the word length by $i$. 
The \emph{differential length of $(\land V,d)$,} denoted $\dl(X)$,
is the
least integer $n$ such that $d_{n-1}$ is non zero. \\
As a minimal model of $X$ is defined up to isomorphism, the differential length is a rational homotopy type invariant of $X$, see \cite[Theorem 1.1]{Kot04}.
\propref{prop:dlnH} establishes a relation between $\dl(X)$ and the existence of an $H(n)$-structure on the rationalisation of $X$.

Finally, recall that \emph{the rational cup-length $\clup_0(X)$ of $X$} is
the maximal length of a nonzero product
in $H^{>0}(X;\mathbb Q)$. In \cite{Kot04}, by using this cup-length and the invariant $\dl(Y)$, Y.~Kotani gives a necessary and sufficient condition for the rational cohomology of
${\mathcal F}_*(X,Y,*)$ to be free commutative when $X$ is a rational formal space and when the dimension of
$X$ is less than the connectivity of $Y$. We show here that a large part of the
Kotani criterium remains valid, without hypothesis of formality and dimension. We prove:

\begin{theorem}\label{thm:rationalcup}
 Let $X$ and $Y$ be nilpotent finite
type CW-complexes, with $X$ finite.
\begin{enumerate}
\item The cohomology
algebra $H^*({\mathcal F}_*(X,Y,*);\mathbb Q)$ is free commutative
if \hfill\break ${\clup}_0(X ) < \dl(Y)$.
 \item If $\dim (X) \leq  \conn (Y)$, then the cohomology
algebra $H^*({\mathcal F}_*(X,Y,*);\mathbb Q)$ is free commutative  if, and only if, $\clup_0(X) <  \dl(Y)$.
 \end{enumerate}
\end{theorem}

As an application, we describe in \thmref{thm:solvable} the Postnikov tower of the rationalisation of ${\mathcal F}_*(X,Y,*)$ where $X$ is a finite nilpotent space and $Y$ a finite type CW-complex whose connectivity is greater than the dimension of $X$. Our description implies the solvability of the rational Pontrjagin algebra of $\Omega({\mathcal F}_*(X,Y,*))$.

\secref{sec:nHspace} contains the topological setting and the proof of \propref{prop:cat}.
The link with rational models is done in \secref{sec:rationalnH}. Our proof of \thmref{thm:rationalcup} uses the Haefliger model for mapping spaces. In order to be self-contained,  we recall briefly Haefliger's construction in \secref{sec:Haefliger}. The proof of \thmref{thm:rationalcup} is contained in \secref{sec:proofrational}. 
Finally,  \secref{sec:Postnikov} is devoted to the description of the Postnikov tower.

In this text, all spaces are supposed of the homotopy type of connected pointed CW-complexes and we  will use cdga for \emph{commutative differential graded algebra.} A \emph{quasi-isomorphism} is a morphism of cdga's which induces an isomorphism in cohomology.

\section{Structure of $H(n)$-space}\label{sec:nHspace}
First we recall the construction of Ganea fibrations, $p_{n}^X\colon G_{n}(X)
\to X$. 
\begin{enumerate}
\item[$\bullet$] Let $F_0(X) \overset {i_0} \to G_0(X) \overset {p_0^X} \to X$
denote the path fibration on $X$, $\Omega X \to PX \to X$.
\item[$\bullet$] Suppose a fibration $\xymatrix@1{F_n(X) \ar[r]^-{i_n} & G_n(X)
\ar[r]^-{p_n^X} & X}$ has been constructed.
We extend $p_n^X$ to a map $q_n \colon  G_n(X)
\cup C(F_n(X))\rightarrow X$, defined on the mapping cone of $i_n$, by setting
$q_n(x)=p_n^X(x)$ for $x \in G_n(X)$ and $q_n([y,t]) = *$
for $[y,t] \in C(F_n(X))$.
\item[$\bullet$] Now convert $q_n$ into a fibration $p_{n+1}^X\colon G_{n+1}(X)
\to X$.
\end{enumerate}

This construction is functorial and the space $G_n(X)$ has the homotopy type of the $n^{\rm th}$-classifying space of Milnor \cite{Mil56b}. We quote also from \cite{Gan67a} that the direct limit $G_\infty(X)$ of the maps $G_n(X)\to G_{n+1}(X)$ has the homotopy type of $X$. As spaces are pointed, one has two canonical applications $\iota_n^l\colon G_n(X)\to G_n(X\times X)$ and $\iota_n^r\colon G_n(X)\to G_n(X\times X)$ obtained from maps $X\to X\times X$ defined respectively by $x\mapsto (x,*)$ and $x\mapsto (*,x)$.

\begin{definition}\label{def:nH}
A space $X$ is \emph{an $H(n)$-space} if there exists a map $\mu_n\colon G_n(X\times X) \to X$ such that $\mu_n\circ \iota_n^l=\mu_n\circ \iota_n^r=p_n^X\colon G_n(X)\to X$.
\end{definition}

Directly from the definition, we see that an $H(\infty)$-space is an $H$-space and that any space is a $H(1)$-space.
Recall also that any co-$H$-space is of LS-category 1. Then, \propref{prop:cat} contains the trivial cases of a co-$H$-space $X$ and of an $H$-space $Y$.

\begin{proof}[Proof of \propref{prop:cat}]
From the hypothesis, we have a section $\sigma\colon X\to G_n(X)$ of the Ganea fibration $p_n^X$  and a map
$\mu_n\colon G_n(Y\times Y)\to Y$ extending the Ganea fibration $p_n^Y$, as in \defref{def:nH}.
If $f$ and $g$ are elements of ${\mathcal F}_*(X,Y,*)$, we set
$$f\bullet g= \mu_n\circ G_n(f\times g)\circ G_n(\Delta_X)\circ \sigma,$$
where $\Delta_X$ denotes the diagonal map of $X$. One checks easily that $f\bullet *\simeq *\bullet f\simeq f$.
\end{proof}

In the rest of this section, we are interested in the existence of $H(n)$-structures on a given space.
For the detection of an $H(n)$-space structure, one may replace the Ganea fibrations $p_n^X$ by any functorial construction of fibrations $\hat{p}_n\colon \hat{G}_n(X)\to X$ such that one has a functorial commutative diagram,
$$\xymatrix{\hat{G}_n(X)\ar@/_.4pc/[rr]\ar[rd]_-{\hat{p}_n}&&\ar@/_.4pc/[ll]G_n(X)\ar[ld]^-{p_n^X}\\
&X.&}$$
 Such maps $\hat{p}_n$ are called fibrations \`a la Ganea in \cite{Sc-Ta97a} and substitutes to Ganea fibrations here. Moreover, as we are interested in product spaces, the following filtration of the space $G_\infty(X)\times G_\infty (Y)$ plays an important role:
$$(G(X)\times G(Y))_n = \cup_{i+j=n} G_i(X)\times
G_j(Y)\,.$$
In \cite{Iwa98}, N. Iwase proved the existence of a commutative diagram
$$\xymatrix{(G(X)\times G(Y))_n\ar@/_.4pc/[rr]\ar[rd]_-{\cup({p}_i^X\times p_j^Y)}&&\ar@/_.4pc/[ll]G_n(X\times Y)\ar[ld]^-{p_n^{X\times Y}}\\
&X\times Y&}$$
 and used it to settle a counter-example to the Ganea conjecture. Therefore, in \defref{def:nH}, we are allowed to replace the Ganea space $G_n(X\times X)$ by $(G(X)\times G(X))_n$. Moreover, if $\hat{p}_n\colon \hat{G}_n(X)\to X$ are substitutes to Ganea fibrations as above, we may also replace $G_n(X\times X)$ by
 $$(\hat{G}(X)\times \hat{G}(Y))_n = \cup_{i+j=n} \hat{G}_i(X)\times
\hat{G}_j(Y)\,.$$
 We will use this possibility in the rational setting.

 In the case $n=2$, we have a cofibration sequence,
$$\xymatrix@1{ \Sigma(G_1(X)\land G_1(X))\ar[r]^-{Wh}& G_1(X)\vee
G_1(X)\ar[r]& G_1(X)\times G_1(X),}$$ coming from the Arkowitz
generalisation of a Whitehead bracket, \cite{Ark62}.  Therefore, the
existence of an $H(2)$-structure on a space $X$ is equivalent to the
triviality of $(p_1^X\vee p_1^X)\circ Wh$. As the loop $\Omega p_1^X$
of the Ganea fibration $p_1^X\colon G_1(X)\to X$ admits a section, we
get the following \emph{necessary condition:}

-- if there is an $H(2)$-structure on $X$, then the homotopy Lie algebra of $X$ is abelian, i.e.\ all Whitehead products vanish.

\begin{example}\label{exam:spheres}
In the case $X$ is a sphere $S^n$, the existence of an $H(2)$ structure on $S^n$ implies $n=1,\,3$ or 7,  \cite{Ada60}. Therefore, only the spheres which are already $H$-spaces endow a structure of  $H(2)$ space. One can also observe that, in general, if a space $X$ is both of category $n$ and an $H(2n)$-space, then it is an $H$-space. The law is given by
$$\xymatrix@1{
X\times X \ar[r]^-{\sigma}&G_{2n}(X\times X)\ar[r]^-{\mu_{2n}}&X},$$
where the existence of the section $\sigma$ to $p_{2n}^{X\times X}$ comes from
$\cat(X\times X)\leq 2\,\cat (X)$.
\end{example}

\begin{example}\label{exam:cp}
If we restrict to spaces whose loop space is a product of spheres or of loop spaces on a sphere, the previous necessary condition becomes a criterion.
For instance, it is proved in  \cite{BJS60} that all Whitehead products are zero in the complex projective 3-space. This implies that $\C P^3$ is an $H(2)$-space.
(Observe that $\C P^3$ is not an $H$-space.)
From  \cite{BJS60}, we know also that the homotopy Lie algebra of $\C P^2$ is not abelian. Therefore $\C P^2$ is not an $H(2)$-space.
\end{example}

The following example shows that we can find $H(n)$-spaces, for any $n>1$.

\begin{example}\label{exam:fiber}
Denote by  $\varphi_r\colon K(\Z,2)\to K(\Z,2r)$ the map corresponding to the class $x^r\in H^{2r}(K(\Z,2);\Z)$, where $x$ is the generator of $H^2(K(\Z,2);\Z)$. Let $E$ be the homotopy fibre of $\varphi_r$. We prove below that \emph{$E$ is an $H(r-1)$-space.}

First we derive, from the homotopy long exact sequence associated to the map $\varphi_r$, that
$\Omega E$ has the homotopy type of $S^1\times K(\Z,2r-2)$. Therefore, the only obstruction to extend
$G_{r-1}(E)\vee G_{r-1}(E)\to E$ to
$(G(E)\times G(E))_{r-1}=\cup_{i+j=r-1}\;G_i(E)\times G_j(E)$
lies in
$${\rm Hom}(H_{2r}((G(E)\times G(E))_{r-1};\Z),\pi_{2r-2}(E)).$$

If $A$ and $B$ are CW-complexes, we denote by $A\sim_{n} B$ the fact that $A$ and $B$ have the same $n$-skeleton. If we look at the Ganea total spaces and fibres, we get:
$$\Sigma\Omega E\sim_{2r} S^2\vee S^{2r-1}\vee S^{2r},\;
F_1(E)=\Omega E*\Omega E\sim_{2r} S^3\vee S^{2r}\vee S^{2r},$$
and more generally, $F_s(E)\sim_{2r} S^{2s+1}$, for any $s$, $2\leq s \leq r-1$. Observe also that $H_{2r}(F_2(E);\Z)\to H_{2r}(G_1(E);\Z)$ is onto. (As we have only spherical classes in this degree, this comes from the homotopy long exact sequence.)

As a conclusion, we have no cell in degree $2r$ in $(G(E)\times G(E))_{r-1}$ and $E$ is an $H(r-1)$-space.
\end{example}

We end this section with a reduction to a more computable invariant than the LS-category. Consider $\rho_n^X\colon X\to G_{[n]}(X)$ the homotopy cofibre of the Ganea fibration~$p_n^X$. Recall that, by definition, $\wcatg (X)\leq n$ if the map $\rho_n^X$ is homotopically trivial. Observe that we always have $\wcatg ( X)\leq \cat( X)$, see \cite[Section~2.6]{CLOT03} for more details on this invariant.

\begin{proposition}\label{prop:wcatg}
Let $X$ be a CW-complex of dimension~$k$ and $Y$ be a CW-complex $(c-1)$-connected with $k\leq c-1$. If $Y$ is an $H(n)$-space such that $\wcatg(X)\leq n$, then ${\mathcal F}_*(X,Y,*)$ is an $H$-space.
\end{proposition}

\proof
Let $f$ and $g$ be elements of ${\mathcal F}_*(X,Y,*)$. Denote by $\tilde{\iota}_n^X\colon \tilde{F}_n(X)\to X$ the homotopy fibre of $\rho_n^X\colon X\to G_{[n]}(X)$. This construction is functorial and the map $(f,g)\colon X\to Y\times Y$ induces a map
$\tilde{F}_n(f,g)\colon \tilde{F}_n(X)\to \tilde{F}_n(Y\times Y)$ such that
$\tilde{\iota}_n^{Y\times Y}\circ \tilde{F}_n(f,g)=(f,g)\circ \tilde{\iota}_n^X$.

By hypothesis, we have a homotopy section
$\tilde{\sigma}\colon X\to \tilde{F}_n(X)$ of $\tilde{\iota}_n^X$. Therefore, one gets a map
$X\to \tilde{F}_n(Y\times Y)$ as
$\tilde{F}_n(f,g)\circ \tilde{\sigma}$.

Recall now that, if $A\to B\to C$ is a cofibration with $A$ $(a-1)$-connected and $C$ $(c-1)$-connected, then  the canonical map $A\to F$ in the homotopy fibre of $B\to C$ is an $(a+c-2)$-equivalence. We apply it in the following situation:
$$\xymatrix{
G_n(Y\times Y)\ar[d]_{j_n^{Y\times Y}}\ar[rr]^{p_n^{Y\times Y}}&&Y\times Y\ar[rr]^{\rho_n^{Y\times Y}}&&G_{[n]}(Y\times Y)\\
\tilde{F}_n(Y\times Y)\ar[urr]_{\tilde{\iota}_n^{Y\times Y}}&&&&\\
}$$
The space $G_n(Y\times Y)$ is $(c-1)$-connected and $G_{[n]}(Y\times Y)$ is $c$-connected. Therefore the map
$j_n^{Y\times Y}$ is $(2c-1)$-connected. From the hypothesis, we get $k\leq c-1 < 2c-1$ and the map $j_n^{Y\times Y}$ induces a bijection
$$\xymatrix@1{[X,G_n(Y\times Y)]\ar[rr]_-{\cong}&&[X,\tilde{F}_n(Y\times Y)].}$$
Denote by $g_n\colon X \to G_n(Y\times Y)$ the unique lifting of $\tilde{F}_n(f,g)\circ \tilde{\sigma}$. The composition $g\bullet f$ is defined as $\mu_n\circ g_n$ where $\mu_n$ is the $H(n)$-structure on $Y$.

If we set $g=*$, then $\tilde{F}_n(f,g)$ is obtained as the composite of $\tilde{F}_n(f)$ with the map 
$\tilde{F}_n(Y)\to \tilde{F}_n(Y\times Y)$
induced by $y\mapsto (y,*)$. As before, one has an isomorphism
$$\xymatrix@1{[X,G_n(Y)]\ar[rr]_-{\cong}&&
[X,\tilde{F}_n(Y)].}$$
A chase in the following diagram shows that $f\bullet *=f$ as expected,
$$\xymatrix{&&G_n(Y)\ar[d]\ar[r]&G_n(Y\times Y)\ar[d]\\
\tilde{F}_n(X)\ar[rr]^{\tilde{F}_n(f)}&&\tilde{F}_n(Y)\ar[d]^{\tilde{\iota}_n^Y}\ar[r]&\tilde{F}_n(Y\times Y)\\
X\ar[u]^-{\tilde{\sigma}}\ar[rr]^f&&Y.&
}\eqno{\raise -80pt\hbox{\qed}}$$

\section{Rational characterisation of $H(n)$-spaces}\label{sec:rationalnH}

Define $m_H(X)$ as the greatest integer $n$ such that $X$ admits an $H(n)$-structure and denote by  $X_0$ the rationalisation of a nilpotent finite type CW-complex $X$. Recall that $\dl(X)$ is the valuation of the differential of the minimal model of $X$, already defined in the introduction. 

\begin{proposition}\label{prop:dlnH}
 Let $X$ be a nilpotent finite type CW-complex of rationalisation $X_0$. Then we have:
 $$m_H(X_0)+1 = \dl(X).$$
\end{proposition}

\begin{proof} Let $(\land V,d)$ be the
minimal model of $X$. Recall from \cite{Fe-Ha82} that a model of the Ganea fibration $p_n^X$ is given by the following composition, 
$$(\land V,d) \to (\land V/\land^{>n}V, \bar{d})\hookrightarrow (\land V/\land^{>n}V, \bar{d})\oplus S,$$
where the first map is the natural projection and the second one the canonical injection together with $S\cdot S
= S\cdot V = 0$ and $d(S) = 0$. As the first map is functorial and the second one admits a left inverse over $(\land V,d)$, we may use the realisation of $(\land V,d) \to (\land V/\land^{>n}V, d)$ as substitute for the Ganea fibration.

Suppose $\dl(X) = r$. We consider the cdga
$(\land V',d') \otimes (\land V'',d'')/I_r$ where $(\land V',d')$ and $(\land V'',d'')$ are copies of $(\land V,d)$ and where $I_r$ is the ideal $I_r =
\oplus_{i+j \geq r} \land^iV'\otimes \land^jV''$. Observe that this cdga has a zero differential and that the morphism
$$\varphi : (\land V,d) \to (\land V',d') \otimes (\land V'',d'')/I_r$$
defined by $\varphi (v) = v'+v''$ satisfies $\varphi(dv)=0$. Therefore $\varphi$  is a morphism of cdga's and its realisation induces an $H(n)$-structure on the rationalisation $X_0$. That shows: $m_H(X_0)+1 \geq \dl(X)$.

Suppose now that $m_H(X_0)+1 > \dl(X)=r$. By hypothesis,
 we have a morphism of cdga's
$$\varphi : (\land V,d) \to (\land V',d') \otimes (\land V'',d'')/I_{r+1}\,.$$ 
By construction, in this quotient, a cocycle of wedge
degree $r$ cannot be a coboundary.
Since the composition of
$\varphi$ with the projection on the two factors is the natural
projection, we have $\varphi(v) -v'-v''\in  \land^+V'\otimes \land^+V''$. Now let $v\in V$, of lowest degree  with $d_r(v) \neq 0$. From $d_r(v) =
\sum_{i_1,i_2,\ldots , i_r}c_{i_1i_2\ldots i_r} v_{i_1}v_{i_2}
\ldots v_{i_r}$, we get
$$\varphi (dv) =
\sum_{i_1,i_2,\ldots , i_r}c_{i_1i_2\ldots i_r}(v_{i_1}'+
v_{i_1}'')\cdot (v_{i_2}'+v_{i_2}'') \cdots (v_{i_r}'+v_{i_r}'')\,.$$
This expression cannot be a coboundary and the equation $d\varphi
(x) = \varphi (dx)$ is impossible. We get a contradiction, therefore one has
$m_H(X_0) + 1 = \dl(X)$.
\end{proof}

 \section{The Haefliger model}\label{sec:Haefliger}

Let $X$ and $Y$ be finite type nilpotent CW-complexes with $X$ of finite dimension.
 Let $(\land V,d)$ be the minimal model of $Y$and $(A,d_A)$ be a finite dimensional model for $X$, which means that $(A,d_A)$ is a finite dimensional cdga equipped with a quasi-isomorphism $\psi $ from the minimal model of $X$ into $(A,d_A)$.
Denote by $A^\vee$ the dual vector space of $A$, graded by
$$(A^\vee)^{-n} = \mbox{Hom} (A^n,\mathbb Q)\,.$$
 We set $A^+ = \oplus_{i=1}^\infty A^i$, and we fix an homogeneous basis $(a_1, \cdots ,a_p)$ of $A^+$. The dual basis $(a^s)_{1\leq s\leq p}$ is a basis of $B = (A^+)^\vee$ defined by
 $\langle a^s;a_t\rangle = \delta_{st}$.
 
  We construct now a morphism of algebras
 $$\varphi : \land V \to A \otimes \land (B\otimes V)$$
by $$\varphi (v) = \sum_{s=1}^p a_s\otimes (a^s\otimes v)\,.$$
In \cite{Hae82} Haefliger proves that there is a unique differential $D$ on $\land (B \otimes V)$ such that $\varphi$ is a morphism of cdga's, i.e.\
$(d_A\otimes D)\circ \varphi = \varphi \circ d$.

In general, the cdga $(\land (B\otimes V), D)$ is not positively graded. Denote by 
$D_0\colon B\otimes V\to B\otimes V$ the linear part of the differential $D$. We define  a cdga $(\land Z,D)$ by constructing $Z$ as the quotient of $B\otimes V$ by  $\oplus_{j\leq 0} (B\otimes V)^j$ and  their image by $D_0$. Haefliger proves:

\begin{theorem}{\rm\cite{Hae82}}\qua The commutative differential graded algebra $(\land Z,D)$ is a model of the mapping space  ${\mathcal F}_*(X,Y,*)$.
\end{theorem}

\section{Proof of \thmref{thm:rationalcup}}\label{sec:proofrational}

\begin{proof}
We start with an explicit description of the Haefliger model, keeping the notation of \secref{sec:Haefliger}.
The cdga $(\land V,d)$ is a minimal model of $Y$ and we choose for $V$ a basis $(v_k)$, indexed by a well-ordered set and satisfying $d(v_k) \in \land (v_r)_{r< k}$ for all $k$.
As homogeneous basis $(a_s)_{1\leq s \leq p}$ of $A$, we choose elements $h_i$, $e_j$ and $b_j$ such that:

-- the elements $h_i$ are cocycles and their classes $[h_i]$  form a linear basis of the reduced cohomology of $A$; 

-- the elements $e_j$ form a linear basis of  a supplement of the vector space of cocycles in $A$, and $b_j = d_A(e_j)$. 

We denote by $h^i$, $e^j$ and $b^j$ the corresponding elements of the basis of $B = (A^+)^\vee$. By developing 
$D_0(\sum_s a_s\otimes (a^s\otimes v))=0$, we get a direct description of the linear part $D_0$ of the differentiel $D$ of the Haefliger model:
$$D_0(b^j\otimes v ) = -(-1)^{\vert b^j\vert} e^j\otimes v \text{ and } D_0(h^i\otimes v ) = 0,\text{ for each }v\in V.$$
A linear basis of the graded vector space $Z$ is therefore given by the elements:
 $$\left\{ \begin{array}{ll} b^j\otimes v_k, & \text{ with } \vert b^j\otimes v_k\vert \geq 1,\\
e^j\otimes v_k, & \text{ with } \vert e^j\otimes v_k\vert \geq 2,\\
h^i\otimes v_k, & \text{ with } \vert h^i\otimes v_k\vert \geq 1.
 \end{array}  \right.$$
 Now, from $\varphi(dv)=(D-D_0)\varphi(v)$ and $d(v) = \sum c_{i_1i_2\cdots  i_r}v_{i_1}v_{i_2}\cdots v_{i_r}$,  we deduce:
 $$\begin{array}{l}
(D-D_0)(a^s\otimes v) =\\
\!\! \pm\sum c_{i_1i_2\cdots i_r}\sum_{a_{i_1}, a_{i_2}\ldots ,a_{i_r}} \langle a^s; a_{i_1}a_{i_2}\cdots a_{i_r}\rangle \, (a_{i_1}\otimes v_{i_1})\cdot (a_{i_2}\otimes v_{i_2})\cdots (a_{i_r}\otimes  v_{i_r})
\end{array}$$
where, as usual, the sign $\pm$ is entirely determined by a strict application of the Koszul rule for a permutation of graded objects.

\medskip

\textsl{Suppose first that $\clup_0(X) < \dl(Y)$.}

We prove, by induction on $k$, that in $(\land Z,D)$ the ideal $I_k$ generated by the elements $$ \left\{
 \begin{array}{ll} 
 b^j\otimes v_s\,, \hspace{3mm} s\leq k,& \text{ with degree at least 1,}\\ e^j\otimes v_s\,, \hspace{3mm}
  s\leq k, & \text{ with degree at least 2,}
  \end{array}
   \right.$$
is a differential ideal and that the elements $h^i \otimes v_s$, with $s\leq k$ and $\vert h^i \otimes v_s\vert \geq 1$, are cocycles in the quotient $((\land Z)/I_k, \bar{D})$.
Note that this ideal is acyclic as shown by the formula given for $D_0$. Therefore the quotient map $\rho\colon (\land Z,D)\to ((\land Z)/I_k, {D})$ is a quasi-isomorphism. The induction will prove that the differential is zero in the quotient, which is the first assertion of \thmref{thm:rationalcup}.

Begin with the induction. 
One has $dv_1=0$ which implies $(D-D_0)(a^s\otimes v_1)=0$. Therefore, we deduce
 $D(b^j\otimes v_1) = -(-1)^{\vert b^j\vert} e^j\otimes v_1$ and $D(h^i\otimes v_1) = 0$. That proves the assertion for $k=1$.

We suppose now that the induction step is true for the integer $k$.
Taking the quotient by the ideal $I_k$ gives a quasi-isomorphism
 $$ \rho : (\land Z,D) \to (\land T,D) :=((\land Z) / I_k, D)\,.$$
As the elements $b^j\otimes v_s$ and $e^j\otimes v_s$, $s\leq k$, have disappeared and as  $\clup_0(X) < \dl(Y)$,  we have $\rho\circ \varphi (dv_{k+1}) = 0$. Therefore
 $D(b^j\otimes v_{k+1}) = -(-1)^{\vert b^j\vert} e^j\otimes v_{k+1}$ and $D(h^i\otimes v_{k+1}) = 0$.
 The induction is thus proved.

\medskip
 \textsl{We consider now the case $\clup_0(X) \geq \dl(Y)$ in the case $\dim (X) \leq  \conn (Y)$.}

We choose first in the lowest possible degree $q$ an element
$y\in V$ that satisfies $ dy = d_ry +\cdots $ with $d_r(y)\neq 0$
and $r\leq \clup_0(X)$. As above we can kill all the elements $e^j\otimes v$ and $b^j\otimes v$ with $\vert v\vert <q$ and keep a quasi-isomorphism $\rho\colon  (\land Z,D)\to (\land T,D):=(\land Z/I_{q-1},D)$.

Next we choose cocycles, $h_1, h_2, \cdots , h_m$, such that the class $[\omega]$, associated to the product $\omega = h_1\cdot h_2 \cdots h_m$, is not zero. We choose
 $m\geq r$ and suppose that $\omega$ is in the highest degree for such a product.
Let
$\omega'$ be an element in $A^\vee$ such that $\langle \omega';\omega\rangle = 1$.
Then, by the Haefliger formula, the differential $D$ is zero in $\land T$ in degrees strictly less than $\vert \omega'\otimes y\vert$. Observe that $\vert \omega'\otimes y\vert\geq 2$ and that the $D_r$ part of the differential $D(\omega'\otimes y)$ is a nonzero sum. This proves that the cohomology is not free.
\end{proof}

\begin{example} Let $X$ be a space with $\clup_0(X) = 1$, which means that all products are zero in the reduced rational cohomology of $X$. Then, for any nilpotent finite type CW-complex $Y$, the rational cohomology $H^*({\mathcal F}_*(X,Y,*);\mathbb Q)$ is a free commutative graded algebra.
For instance, this is the case for the (non-formal) space $X = S^3\vee S^3 \cup_\omega e^8$, where the cell $e^8$ is attached along a sum of triple Whitehead products.
\end{example}

\begin{example} When the dimension of $X$ is greater than the connectivity of $Y$,
the degrees of the elements have some importance. The cohomology can be commutative free  even if $\clup_0(X)\geq \dl(Y)$.
For instance, consider $X = S^5 \times S^{11}$ and $Y = S^8$. One has
$\clup_0(X)=\dl(Y)=2$ and the function space ${\mathcal F}_*(X,Y,*)$ is a rational $H$-space with the rational homotopy type of $K(\mathbb Q,3)\times K(\mathbb Q,4)\times K(\mathbb Q, 10)$, as a direct computation with the Haefliger model shows.
\end{example} 

\section{Rationalisation of ${\mathcal F}_*(X,Y,*)$ when 
$\dim(X)\leq$\hfill\break$ \conn(Y)$} \label{sec:Postnikov}

Let $X$ be a finite nilpotent space with rational LS-category equal to $m-1$
and let $Y$ be a finite type nilpotent CW-complex whose connectivity $c $ is greater than the dimension of $X$.
We set $r = \dl(Y)$ and denote by $s$ the maximal integer
such that $m/r^s \geq 1$, i.e.\ $s$ is the integral part of $\log_rm$.

\begin{theorem}\label{thm:solvable}
There is a sequence of rational fibrations $K_k\to F_k\to F_{k-1}$, for $k = 1, \ldots, s$, with $F_0 = *$,
$F_s$ is the rationalisation of ${\mathcal F}_*(X,Y,*)$ and each space $K_k$ 
is a product of Eilenberg-MacLane spaces.
In particular, the rational loop space homology of ${\mathcal F}_*(X,Y,*)$ is solvable with solvable
index less than or equal to $s$.
\end{theorem}

\begin{proof}
 By a result of Cornea \cite{Cor94a}, the space $X$ admits a finite dimensional model $A$ such that $m$ is the maximal length of a nonzero product of elements of positive degree. We denote by $(\land V,d)$ the minimal model of $Y$.

We consider the ideals $I_k = A^{> m/r^k}$, and the short exact sequences of cdga's
$$I_k/I_{k-1} \to A/I_{k-1} \to A/I_k\,.$$
These short exact sequences realise into cofibrations $T_k \to T_{k-1} \to Z_k$ and the sequences
$$(\land ((A^+/I_k)^\vee \otimes V),D) \to (\land ((A^+/I_{k-1})^\vee \otimes V),D) \to
(\land ((I_k/I_{k-1})^\vee \otimes V),D)$$
are relative Sullivan models for the fibrations
$${\mathcal F}_*(Z_k, Y,*) \to {\mathcal F}_*(T_{k-1}, Y,*)\to {\mathcal F}_*(T_k, Y,*).$$
Now since the cup length of the space $Z_k$ is strictly less than $r$, the function spaces
${\mathcal F}_*(Z_k, Y,*) $ are rational $H$-spaces, and this proves \thmref{thm:solvable}.
\end{proof}

\Addresses\recd

\end{document}